\documentclass[12pt]{amsart}
\usepackage{amsmath}
\usepackage{amssymb}
\usepackage{latexsym}

\newenvironment{proof*}{\vskip 2mm\noindent {}}{\hfill $\Box$ \vskip 2mm}
\newtheorem{theorem}{Theorem}

\newtheorem{definition}[theorem]{Definition}

\newcommand{\C}{{\mathbb{C}}}

\newcommand{\Ker}{\operatorname{Ker}}

\title[Automorphisms of the spectral ball]{A local form for the automorphisms of the spectral unit ball}
\author{Pascal J. Thomas}

\address{Institut de Math\'ematiques de Toulouse, UMR CNRS 5219\\
Universit\'e Paul Sabatier, 118 Route de Narbonne\\ F-31062 Toulouse
Cedex, France} \email{pthomas@cict.fr}



\begin{document}

\begin{thanks}{This paper was made possible in part by a grant from the French Ministry of
Foreign Affairs in the framework of the ECO Net programme, file 10291 SL, coordinated by Ahmed Zeriahi.}
\end{thanks}

\begin{abstract} If $F$ is an automorphism of $\Omega_n,$ the $n^2$-dimensional
spectral unit ball, we show that, in a neighborhood of any cyclic matrix of $\Omega_n,$ the map $F$ can be written as conjugation by a holomorphically varying non singular matrix. This provides a shorter proof of a theorem of J. Rostand, with a slightly stronger result.
\end{abstract}

\maketitle

\section{Background}

Let $\mathcal M_n$ be the set of all $n\times n$ complex matrices.
For $A\in\mathcal M_n$ denote by $sp(A)$ the spectrum of $A.$ 
The spectral ball $\Omega_n$ is the set
$$\Omega_n:=\{A\in\mathcal M_n: \forall \lambda\in sp(A), |\lambda|<1\}.$$

Let $F$ be an automorphism of $\Omega_n$, that is to say, a biholomorphic map
of the spectral ball into itself. Ransford and White \cite{RW} proved that, by composing
with a natural lifting of a M\"obius map of the disk, one could reduce oneself to the case
where $F(0)=0$, and that in that case the linear map $F'(0)$ was a linear automorphism
of $\Omega_n,$ so that by composing with its inverse, one is reduced to the case $F(0)=0$,
$F'(0)=I$ (the identity map). We then say that the automorphism if \emph{normalized}.
Ransford and White \cite{RW} proved that such automorphisms preserve the spectrum of matrices.

We say that two matrices $X, Y$ are \emph{conjugate} if there exists $Q \in \mathcal M_n^{-1}$
such that $X = Q^{-1} Y Q$. 

Baribeau and Ransford \cite{BRa} (see also \cite{BRo} for a more elementary proof) proved that 
every spectrum-preserving $\mathcal C^1$-diffeomorphism of an open subset of $\mathcal M_n$, and thus every
normalized automorphism of the spectral ball is a pointwise conjugation: 
\begin{equation}
\label{conj}
F(X) = Q(X)^{-1} X Q(X).
\end{equation}
Rostand's contribution \cite{Ro} was to show that $Q(X)$ could be chosen locally holomorphically
in a neighborhood of every $X$ admitting $n$ distinct eigenvalues. 

We will give a short proof of a slightly stronger result: the exceptional set of matrices where
the local holomorphic choice cannot be guaranteed will be of complex
codimension $2$ instead of $1$.

The motivation for this result was a conjecture formulated in \cite{RW} about the 
automorphisms of the spectral ball, which reduces to asking whether any normalized automorphisms
can be written in the form \eqref{conj}, where $Q$ would be globally homorphic on $\Omega_n$,
and depend only on the conjugacy class of $X$. Notice that a recent result of Zwonek \cite{Zwo}
shows that any proper map of the spectral ball to itself is actually an automorphism of it, so that 
the proof of the Ransford-White conjecture would yield a description of all the proper maps
of the spectral ball into itself. 

I wish to thank Nikolai Nikolov, who told me about this circle of ideas. Without the fruitful discussions I had with him on this topic, this paper wouldn't have been written. 

\section{Statement}

\begin{definition}
We say that a matrix $M$ is \emph{cyclic} (or \emph{non-derogatory}) if there exists a cyclic vector for $M$,
i.e. $v \in \mathbb C^n$ such that $(v, Mv, \dots, M^k v, \dots)$ spans $\mathbb C^n$,
which is equivalent to the fact that $(v, Mv, \dots, M^{n-1} v)$ is a basis of $\mathbb C^n$.
\end{definition}

Many equivalent definitions of this notion can be found, for instance in \cite{HJ}
and \cite{HJTop}, or \cite[Proposition 3]{NTZ}.
We point one out: $M$ is cyclic if and only if for any $\lambda \in \mathbb C$,
$\dim \Ker(M-\lambda I_n) \le 1$. In particular, any matrix with $n$
distinct eigenvalues is cyclic,  and for any given spectrum $\lambda_1, \dots, \lambda_n$,
the set of non-cyclic matrices with that spectrum is the algebraic set
$$
\{ M : \exists j : \dim \Ker(M-\lambda_j I_n) \ge 2 \}.
$$
Hence the set of non-cyclic matrices is of codimension $1$ in the set of matrices
which admit at least one multiple eigenvalue, itself of codimension $1$ in $\mathcal M_n$.

\begin{theorem}
Let $F$ be a spectrum-preserving holomorphic map of $\Omega_n$. 
Let $X_0 \in \Omega_n$ be a cyclic matrix. Then there exists a neighborhood $\mathcal V_{X_0}$ of $X_0$
and a map $Q$ holomorphic from $\mathcal V_{X_0}$ to $\mathcal M_n^{-1}$ such that for 
any $X \in \mathcal V_{X_0}$, $F(X) = Q(X)^{-1} X Q(X)$.
\end{theorem}

\section{Proof}

We fix some notation.
For $A\in\mathcal M_n,$ let
$$\sigma_j(A):=\sigma_j(\lambda_1,\dots,\lambda_n):=
(-1)^j \sum_{1\le k_1<\dots<k_j\le n}\lambda_{k_1}\dots\lambda_{k_j}$$ and
$\lambda_1,\dots,\lambda_n$ are the (possibly equal) eigenvalues of $A.$ Those are polynomials in the 
coefficients of $A$. 
Put $\sigma:=(\sigma_1,\dots,\sigma_n):\mathcal M_n\to \mathbb C^n.$

Conversely, given $a:= (a_1, \ldots, a_{n}) \in \C^n$, the associated \emph{companion matrix} 
$\mathcal C_a$ is
$$
\left( 
\begin{array}{ccccc}
0&&&&-a_n\\
1&0&&&\vdots\\
&1&\ddots&&\vdots\\
&&\ddots&0&-a_{2}\\
&&&1&-a_{1}
\end{array}
\right) .
$$
The companion matrix associated to a matrix $A$ is $\mathcal C_{\sigma(A)}$. They have the
same characteristic polynomial, or equivalently $\sigma(\mathcal C_{\sigma(A)}) = \sigma(A)$.

Now given a matrix $X_0$ as in the Theorem, and a vector $v_0$ cyclic for $X_0$, let
$$
\mathcal U_{X_0} := \{ M \in \mathcal M_n : det ( v_0, M v_0, \dots, M^{n-1} v_0 ) \neq 0 \}.
$$
This is a neighborhood of $X_0$. Let $P_{v_0}(M)$ be the matrix with columns $( v_0, M v_0, \dots, M^{n-1} v_0 )$;
this depends polynomially on the entries of $M$, and is invertible. One can see that
for $X \in \mathcal U_{X_0}$,
$$
P_{v_0}(X)^{-1} X P_{v_0}(X) = \mathcal C_{\sigma(X)}
$$
(the $n-1$st columns columns coincide, and they have the same characteristic polynomial).

By the Baribeau-Ransford theorem \cite{BRa} $F(X_0)$ is conjugate to $X_0$, therefore
cyclic.
So there is a neighborhood $\mathcal U_{F(X_0)} $ where the relation 
\newline
$P_{w_0}(Y)^{-1} Y P_{w_0}(Y) = \mathcal C_{\sigma(Y)}$ holds. Take $\mathcal V_{X_0}
\subset \mathcal U_{X_0}$ small enough so that $F(\mathcal V_{X_0}) \subset \mathcal U_{F(X_0)} $.
For any $X \in \mathcal V_{X_0}$, using the fact that $F$ is spectrum-preserving, 
\begin{multline*}
F(X) =  P_{w_0}(F(X)) \mathcal C_{\sigma(F(X))} P_{w_0}(F(X))^{-1}\\
= P_{w_0}(F(X)) \mathcal C_{\sigma(X)} P_{w_0}(F(X))^{-1} \\
= P_{w_0}(F(X))P_{v_0}(X)^{-1} X P_{v_0}(X) P_{w_0}(F(X))^{-1},
\end{multline*}
so that we have the theorem with $Q(X) = P_{v_0}(X) P_{w_0}(F(X))^{-1}$.
$\square$

Notice that there is no hope to make a global holomorphic choice of $v_0$ on the 
whole open of cyclic matrices. Indeed, since the complement of this is of codimension $2$,
we could then extend it to the whole of $\Omega_n$ by Hartog's phenomenon, but it would mean that all matrices
are cyclic, which is obviously false.


\begin{thebibliography}{}

\bibitem{BRa} L. Baribeau and T. J. Ransford, {\it Non-linear spectrum-preserving maps},
Bull. London Math. Soc., 32 (2000), 8--14.

\bibitem{BRo} L. Baribeau and S. Roy, {\it Caract\'erisation spectrale de la forme de Jordan},
Linear ALgebra Appl. 320 (2000), 183--191. 

\bibitem{HJ} R. A. Horn, C. R. Johnson, {\it Matrix Analysis}, Cambridge University
Press, Cambridge, New York, Melbourne, 1985.

\bibitem{HJTop} R. A. Horn, C. R. Johnson, {\it Topics in Matrix Analysis}, Cambridge University
Press, Cambridge, New York, Melbourne, 1991.

\bibitem{NTZ} N. Nikolov, P. J. Thomas, W. Zwonek,
{\it Discontinuity of the Lempert function and the Kobayashi--Royden
metric of the spectral ball}, Preprint, 2007
(arXiv:math.CV/0704.2470).


\bibitem{RW} T. J. Ransford and M. C. White, {\it Holomorphic self-maps of the spectral unit ball},
Bull. London Math. Soc., 23 (1991), 256--262.

\bibitem{Ro} J. Rostand, {\it On the automorphisms of the spectral unit ball}, Studia Math. 155 (3)
(2003), 207--230.

\bibitem{Zwo} W. Zwonek, {\it  Proper holomorphic mappings of the spectral unit ball},
preprint, 2007
(arXiv:math.CV/0704.0614).

\end{thebibliography}
\end{document}